\input amstex
\input amstex
\input amsppt1
\documentstyle{amsppt}
\nologo
\redefine\R{\Bbb{R}}
\topmatter
\title about the injectivity radius and the Ricci tensor of a complete Riemannian Manifold\endtitle
\author S\'ergio L. Silva\endauthor
\address Departamento de Estruturas Matem\'aticas, IME, Universidade do Estado do Rio de Janeiro, 20550-013, Rio de Janeiro, Brazil.\endaddress
\email sergiol\@ime.uerj.br\endemail
\subjclass\nofrills 2000 {\it Mathematics Subject Classification.} Primary ; Secondary \endsubjclass
\abstract In this paper we obtain a simple upper bound for the infimum of the Ricci curvatures of a complete Riemannian manifold $M^n$ with nonzero injectivity radius $i(M)$ depending only on of the $i(M)$.  In case of rigidity the Riemannian manifold must be an $n$-dimensional Euclidean sphere(Euclidean space) conform the injectivity radius be finite(infinite). Furthermore with the additional assumption that the second derivative of the Ricci tensor of $M^n$ is null we prove that the same upper bound for the infimum of the Ricci curvatures holds for the supremum of the Ricci curvatures and $M^n$ has, in fact, parallel Ricci tensor. \endabstract
\endtopmatter

\document

\rightheadtext{injectivity radius and the Ricci tensor}

\head \bf 1. Introduction \endhead

Let $M^n$ be a differentiable $n$-dimensional complete Riemannian manifold with metric $\langle\,,\,\rangle$ and Levi-Civita connection $\nabla$. The {\it curvature tensor} $R$ of $M^n$ is defined by
$$
R(X,Y)Z=\nabla_X\nabla_Y Z-\nabla_Y\nabla_X Z-\nabla_{[X,Y]}Z
$$
for all differentiable vector fields $X,\,Y$ and $Z$ on $M^n$.

We denote by $T_pM$ the tangent space of $M^n$ at $p$, $exp_p$ the exponential map of $M^n$ at $p$, $d$ the intrinsic distance function on $M^n$ an by $C(p)$ the {\it cut locus of $p$}. The {\it injectivity radius of $M^n$} is given by
$$
i(M)=\underset {p\in M} \to \inf\,d(p,C(p)).
$$ 

The {\it Ricci Tensor} of $M^n$, denoted by $Ric$, is the differentiable symmetric $2$-form
$$
Ric_p(v,w)=\frac{1}{n-1}\sum_{\imath=1}^n \left\langle R_p(v,e_\imath)e_\imath,w\right\rangle
$$
for all $p\in T_pM$ and $v,\,w\in T_pM$, being $e_1,\ldots, e_n$ an orthonormal basis of $T_pM$. Given $p\in M^n$ and a unitary vector $v$ tangent to $M^n$ at $p$, the number $Ric_p(v,v)$ is called {\it the Ricci curvature of $M^n$ at $p$ in the direction $v$}.

For a differentiable $k$-form $\omega$ on $M^n$ its covariant derivative is the $(k+1)$-form $\nabla \omega$ given by
{\eightpoint$$
\nabla \omega\left(X_1,X_2,\ldots,X_{k+1}\right)=X_1\left(\omega\left(X_2,\ldots,X_{k+1}\right)\right)-\sum_{\imath=2}^{k+1}\omega\left(X_2,\ldots,\nabla_{X_1}X_{\imath},\ldots, X_{k+1}\right)
$$}
for all vector fields $X_1,\,X_2,\ldots, X_{k+1}$ on $M^n$. For $m\geq 2$, we put $\nabla^m\omega=\nabla\left(\nabla^{m-1}\omega\right)$. We say that $\omega$ is {\it parallel} on $M^n$ when $\nabla \omega\equiv 0$.

Let
 $T_1M$ denote the unitary tangent bundle of $M^n$, that is, 
$$T_1M=\left\{\,(q,w)\;:\; q\in M^n,\;w\in T_qM,\,||v||=1\,\right\}.$$
Our first result is a simple upper bound for the infimum of the Ricci curvature of a complete Riemannian manifold with nonzero injectivity radius as a function on $T_1M$. If $i(M)$ is infinite follows immediately from Myers' and Bonnet's Theorems that $\underset{\eightpoint (q,w)\in T_1M}\to{\inf}Ric_q(w,w)\leq 0$. Note that if equality holds then $M^n$ is isometric to the Euclidean space $\Bbb R^n$ by a splitting theorem(see Theorem 2 in \cite{C-G}) since any geodesic is a line. When $i(M)$ is finite we obtain the following result

\par\proclaim{Theorem 1.1} Let $M^n$ be a complete Riemannian manifold with $0<i(M)<\infty$. Then, 
$$
\underset{\eightpoint (q,w)\in T_1M}\to{\inf}Ric_q(w,w)\leq \left[\frac{\pi}{i(M)}\right]^2.
$$
If equality holds then $M^n$ is isometric to the sphere of diameter $i(M)$. 
\endproclaim

Our next result establishes that if a complete Riemannian manifold satisfies $\nabla^2Ric\equiv 0$  then $M^n$ has, in fact, parallel Ricci tensor and the same upper bound for the infimum of the Ricci curvatures also holds for the supremum of the Ricci curvature. More precisely, we have the following result

\par\proclaim{Theorem 1.2} Let $M^n$ be a complete Riemannian manifold with $\nabla^2Ric\equiv 0$.  Then the following hold:

\bigskip

(i) If $i(M)\neq 0$ then the Ricci tensor of $M^n$ is parallel. Furthermore, if $i(M)=\infty$ then the Ricci curvature  is non-positive on $M^n $ and if $0<i(M)<\infty$ then $\underset{\eightpoint (q,w)\in T_1M}\to{\sup}Ric_q(w,w)\leq \left[\frac{\pi}{i(M)}\right]^2$.

\bigskip

(ii) If $n=2$ ($M$ is a surface) then $M^2$ has constant Gaussian curvature.  
\endproclaim

\head \bf 2. Proofs of Theorem 1.1 and Theorem 1.2\endhead

\demo{Proof of Theorem 1.1} Given $(p,v)\in T_1M$ and a real number $t_o$ such that $0<t_o<d(p,C(p))$, let $\gamma$ be the normalized geodesic satisfying $\gamma(0)=p$ and $\gamma^\prime(0)=v$. The geodesic $\gamma$ is minimizing from $0$ to $t_o$. Consider an orthonormal basis $e_1=v$, $e_2$, $\ldots$, $e_n$ of $T_pM$ and the variations $f_i(t,s)=exp_{\gamma(t)}s\,V_i(t)$, $t\in [0,t_o]$, $s\in\R$, where  $V_i(t)=e_i(t)\sin\left(\frac{\pi\,t}{t_o}\right)$ with $e_i(t)$ being the parallel translate of $e_i$ along $\gamma$. For each $i$ and $s$, the curve $\theta_s^i(t)=f_i(t,s)$, $t\in [0,t_o]$, is a curve joining $p$ to $\gamma\left(t_o\right)$ with $\theta_0^i=\gamma|_{[0,t_o]}$. Since $\theta_0^i$ is minimizing we can write
$$
0\leq E_i^{\prime\prime}(0)\;\;\text{for all}\;i=2,\ldots,n. \tag2.1
$$
Above $E_i(s)=\int_0^{t_o}\left|\frac{d\theta_s^i}{dt}(t)\right|^2dt$ is the energy of $\theta_s^i$. The second derivative of $E_i$ is given by
$$\align
\frac{1}{2}E_i^{\prime\prime}(s)&=-\int_0^{t_o}\left\langle R\left(\frac{\partial f_i}{\partial t},\frac{\partial f_i}{\partial s}\right)\frac{\partial f_i}{\partial s},\frac{\partial f_i}{\partial t}\right\rangle dt\tag2.2\\
&+\int_0^{t_o}\left[\left\langle\frac{D}{dt}\frac{D}{ds}\frac{\partial f_i}{\partial s},\frac{\partial f_i}{\partial t}\right\rangle+\left|\frac{D}{dt}\frac{\partial f_i}{\partial s}\right|^2\right] dt.
\endalign$$
We have
$\frac{\partial f_i}{\partial t}(t,0)=\gamma^{\prime}(t)$, $\frac{\partial f_i}{\partial s}(t,0)=V_i(t)$ and $\frac{D}{ds}\frac{\partial f_i}{\partial s}(t,s)=0$. Consequently,
$$\align
\frac{1}{2}E_i^{\prime\prime}(0)&=-\int_0^{t_o}\sin^2\left(\frac{\pi t}{t_o}\right)\left\langle R\left(\gamma^{\prime}(t),e_i(t)\right)e_i(t),\gamma^{\prime}(t)\right\rangle dt\tag2.3\\
&+\int_0^{t_o}\left(\frac{\pi}{t_o}\right)^2\cos^2\left(\frac{\pi t}{t_o}\right)dt.
\endalign$$
and
$$\align
\frac{1}{2}\sum_{\imath =2}^nE_i^{\prime\prime}(0)&=-(n-1)\int_0^{t_o}\sin^2\left(\frac{\pi t}{t_o}\right)Ric_{\gamma(t)}\left(\gamma^{\prime}(t),\gamma^{\prime}(t)\right) dt\tag2.4\\
&+(n-1)\int_0^{t_o}\left(\frac{\pi}{t_o}\right)^2\cos^2\left(\frac{\pi t}{t_o}\right)dt.
\endalign$$
Due to \thetag{2.1} we can write
$$\align
0&\leq -\int_0^{t_o}\sin^2\left(\frac{\pi t}{t_o}\right)Ric_{\gamma(t)}\left(\gamma^{\prime}(t),\gamma^{\prime}(t)\right) dt+\int_0^{t_o}\left(\frac{\pi}{t_o}\right)^2\cos^2\left(\frac{\pi t}{t_o}\right)dt\\
&\leq -\underset{\eightpoint (q,w)\in T_1M}\to{\inf}Ric_q(w,w)\int_0^{t_o}\sin^2\left(\frac{\pi t}{t_o}\right) dt+\int_0^{t_o}\left(\frac{\pi}{t_o}\right)^2\cos^2\left(\frac{\pi t}{t_o}\right)dt\\
&=\frac{t_o}{2}\left[\left(\frac{\pi}{t_o}\right)^2-\underset{\eightpoint (q,w)\in T_1M}\to{\inf}Ric_q(w,w)\right].
\endalign$$
The above inequality implies that $\underset{\eightpoint (q,w)\in T_1M}\to{\inf}Ric_q(w,w)\leq \left(\frac{\pi}{t_o}\right)^2$ for all $t_o$ such that $0<t_o<d(p,C(p))$. If $d(p,C(p))$ is infinite for some $p$ then $\underset{\eightpoint (q,w)\in T_1M}\to{\inf}Ric_q(w,w)\leq 0$ and inequality in Theorem 1.1 is immediate. So we can suppose $d(p,C(p))$ finite for all $p\in M$. In this case, 
$$
\underset{\eightpoint (q,w)\in T_1M}\to{\inf}Ric_q(w,w)\leq \left(\frac{\pi}{d(p,C(p))}\right)^2.
$$ 
Since $i(M)=\underset{p\in M}\to{\inf}d(p,C(p))$ and $i(M)$ is finite and nonzero the inequality in Theorem 1.1 follows. If equality holds then $Ric_q(w,w)\geq \left(\frac{\pi}{i(M)}\right)^2$ for all $(q,w)\in T_1M$. So by Myers' and Bonnet's Theorems $M^n$ is compact and has diameter less than or equal to $i(M)$. Consequently, the diameter of $M^n$ is $i(M)$. Now that $M^n$ is isometric to the sphere of diameter $i(M)$ follows from Theorem 3.1 in \cite{SC}(see also \cite{KS}).\qed
\enddemo 

Before proving Theorem 1.2 we prove the following proposition

\par\proclaim{Proposition 2.1} Let $M^n$ be a complete Riemannian manifold with $\nabla^2Ric\equiv 0$. Then for all $p\in M^n$, all unitary vector $v$ in $T_pM$ and all $t_o$ such that $0<t_o<d(p,C(p))$ we have $Ric_p(v,v)\leq \left(\frac{\pi}{t_o}\right)^2$. In particular, if $d(p,C(p))$ is finite then $Ric_p(v,v)\leq \left(\frac{\pi}{d(p,C(p))}\right)^2$ for all unitary vector $v$ in $T_pM$.\endproclaim

\demo{Proof} Since $\nabla^2Ric \equiv 0$, for the normalized geodesic $\gamma$ such that $\gamma(0)=p$ and $\gamma^\prime(0)=v$, we have 
$$
Ric_{\gamma(t)}\left(\gamma^\prime(t),\gamma^\prime(t)\right)=a_\gamma\,t+b_\gamma,\;t\in \R,\tag2.5
$$
with $a_\gamma=\nabla Ric_p(v,v,v)$ and $b_\gamma=Ric_p(v,v)$. As in proof of Theorem 1.1 if we consider the variations $f_i(t,s)$, $i=2,\ldots,n$, using \thetag{2.1}, \thetag{2.4} and \thetag{2.5}, we deduce that
$$
0\leq-\int_0^{t_o}\sin^2\left(\frac{\pi t}{t_o}\right)(a_\gamma\,t+b_\gamma) dt+\int_0^{t_o}\left(\frac{\pi}{t_o}\right)^2\cos^2\left(\frac{\pi t}{t_o}\right)dt.\tag2.6
$$
Thus,
$$
0\leq -\frac{a_\gamma}{4}t_o^2-\frac{b_\gamma\,t_o}{2}+\frac{1}{2}\frac{\pi^2}{t_o},\tag2.7
$$
that is,
$$
0\leq -a_\gamma\,t_o^2-2b_\gamma\,t_o+2\frac{\pi^2}{t_o}.\tag2.8
$$ 
If we use the geodesic $\beta$ such that $\beta(0)=p$ and $\beta^\prime(0)=-v$, proceeding analogous to the above, give us
$$
0\leq a_\gamma\,t_o^2-2b_\gamma\,t_o+2\frac{\pi^2}{t_o}.\tag2.9
$$
As a consequence of \thetag{2.8} and \thetag{2.9} we have
$$
0\leq -4b_\gamma\,t_o+4\frac{\pi^2}{t_o}.
$$ 
and
$$
Ric_p(v,v)=b_\gamma\leq \left(\frac{\pi}{t_o}\right)^2.\tag2.10
$$
\qed
\enddemo

\par\proclaim{Corollary 2.2} Let $M^n$ be a complete Riemannian manifold with $\nabla^2Ric\equiv 0$. If $\gamma$ is a normalized ray satisfying $\gamma(0)=p$ and $\gamma^\prime(0)=v$ then $a_\gamma=\nabla Ric_p(v,v,v)\leq 0$. Case $a_\gamma=0$,  we must have $b_\gamma=Ric_p(v,v)\leq 0$. Equivalently, $Ric_\gamma\left(\gamma^\prime,\gamma^\prime\right)$ is either decreasing or constant non-positive when $\gamma$ is a ray. \endproclaim

\demo{Proof} Since $\gamma$ is a ray it is minimizing from $0$ to $t_o$ for all $t_o>0$. Then the inequality \thetag{2.8} holds for all $t_o>0$. Consequently, $a_\gamma=\nabla Ric_p(v,v,v)\leq 0$. If $a_\gamma<0$ then $Ric_{\gamma(t)}\left(\gamma^\prime(t),\gamma^\prime(t)\right)=a_\gamma t+b_\gamma,\;t\in \R$, is decreasing. Case $a_\gamma=0$, we have $Ric_{\gamma(t)}\left(\gamma^\prime(t),\gamma^\prime(t)\right)=b_\gamma$ for all $t\in\R$ and is clear from \thetag{2.8} that $b_\gamma\leq 0$.\qed
\enddemo

\par\proclaim{Corollary 2.3} Let $M^n$ be a complete Riemannian manifold with $\nabla^2Ric\equiv 0$. If $\gamma$ is a normalized line satisfying $\gamma(0)=p$ and $\gamma^\prime(0)=v$ then $\nabla Ric_p(v,v,v)=0$ and $Ric_p(v,v)\leq 0$. Equivalently, $Ric_\gamma\left(\gamma^\prime,\gamma^\prime\right)$ is constant non-positive when $\gamma$ is a line. \endproclaim

\demo{Proof} Since $\gamma$ is a line the inequalities \thetag{2.8} and \thetag{2.9} hold for all $t_o>0$. Thus, $a_\gamma=\nabla Ric_p(v,v,v)=0$ and $Ric_{\gamma(t)}\left(\gamma^\prime(t),\gamma^\prime(t)\right)=b_\gamma$ for all $t\in \R$ with $b_\gamma=Ric_p(v,v)\leq 0$.\qed
\enddemo

\par\proclaim{Corollary 2.4} Let $M^n$ be a complete Riemannian manifold with $\nabla^2Ric\equiv 0$. If $d(p,C(p))=\infty$ then $Ric_p(v,v)\leq 0$ for all unitary vector $v\in T_pM$.\endproclaim

\demo{Proof} Since $d(p,C(p))=\infty$, by Proposition 2.1, the inequality $Ric_p(v,v)\leq \left(\frac{\pi}{t_o}\right)^2$ holds for all unitary vector $v\in T_pM$ and all $t_o>0$. Consequently, $Ric_p(v,v)\leq 0$ for all unitary vector $v\in T_pM$.\qed
\enddemo

\bigskip

\demo{Proof of Theorem 1.2} 

\bigskip

(i) If $i(M)=\infty$ then $d(p,C(p))=\infty$ for all $p\in M^n$ and the Corollary 2.4 gives $Ric_p(v,v)\leq 0$ for all $p\in M^n$ and all unitary vector $v\in T_pM$.
If $0<i(M)<\infty$ then for all $p$ such that $d(p,C(p))<\infty$ it holds that $Ric_p(v,v)\leq \left[\frac{\pi}{d(p,C(p))}\right]^2$ for all unitary vector $v\in T_pM$ due to Proposition 2.1. Consequently,  
$Ric_p(v,v)\leq \left[\frac{\pi}{i(M)}\right]^2$ for all $p\in M^n$ and all unitary vector $v\in T_pM$. In any case there exists a nonnegative constant $\rho$ such that $Ric_p(v,v)\leq \rho$ for all $p\in M^n$ and all unitary vector $v\in T_pM$. Now given $p\in M^n$ and unitary vectors $v,\,w\in T_pM$ we consider the geodesic $\gamma$ such that $\gamma(0)=p$ and $\gamma^\prime(0)=v$. If $w(t)$ is the parallel translate of $w$ along $\gamma$, we can write
$Ric_{\gamma(t)}(w(t),w(t))\leq\rho$ for all $t\in \R$. Since $\nabla^2Ric\equiv 0$ we have $Ric_{\gamma(t)}(w(t),w(t))=\nabla Ric_p(v,w,w)\,t+Ric_p(w,w)$. Consequently, $\nabla Ric_p(v,w,w)=0$ for all $p\in M^n$ and all unitary vectors $v,\,w\in T_pM$. Since $\nabla Ric_p$ is a symmetric tensor in the last two coordinates we deduce that $\nabla Ric_p\equiv 0$. Being $p$ arbitrary, we conclude that $M^n$ has parallel Ricci tensor.\qed
  
\bigskip

(ii) If $dim\,M=2$ then the Ricci curvature $Ric_p(v,v)$, for all $p\in M^n$ and all unitary vector $v\in T_pM$, is the Gaussian curvature at $p$ denoted by $K(p)$. So for a normalized geodesic $\gamma$ we have $K_\gamma(t)=K(\gamma(t))=a_\gamma\,t+b_\gamma$, $t\in \R$, where $a_\gamma$ and $b_\gamma$ are real constants depending on $\gamma$. We affirm that $K_\gamma$ is constant for any geodesic $\gamma$. In fact, suppose on the contrary that $a_\gamma\neq 0$ for some geodesic $\gamma$. Changing $\gamma$ by $\beta$ such that $\beta(t)=\gamma(-t)$ if necessary we can assume that $a_\gamma>0$, that is, $K_\gamma$ is increase. Observe that the geodesic $\gamma$ is not contained in any compact set of $M^n$ since $K_\gamma$ is unbounded. Consider $t_o>0$ such that $K_\gamma(t_o)>0$. Since $dK_{\gamma(t_o)}\left(\gamma^\prime(t_o)\right)=a_\gamma>0$ and $dK_{\gamma(t_o)}\left(-\gamma^\prime(t_o)\right)=-a_\gamma<0$, there exists an unitary vector $w\in T_{\gamma(t_o)}M$ such that $dK_{\gamma(t_o)}(w)=0$. Let $\theta$ be the geodesic with $\theta(0)=\gamma(t_o)$ and $\theta^\prime(0)=w$. It holds that $K_\theta(t)=K(\theta(t))=a_\theta\,t+b_\theta$ with $a_\theta=K_\theta^\prime(0)=dK_{\gamma(t_o)}(w)=0$. Consequently, $K_\theta$ is constant. Take a sequence of positive real numbers $t_n$ such that $t_n\to \infty$ and that $p_n=\gamma\left(t_o+t_n\right)$ be divergent. Fixed a point $q=\theta(s)$, consider a normalized minimal geodesic $\gamma_n$ from $q$ to $p_n$. Observe that $b_\theta=K_\theta(s)=K(q)=K_\theta(0)=K(\gamma(t_o))=K_\gamma(t_o)>0$ because $K$ is constant over $\theta$. Supposing $\gamma_n(t)=exp_qtu_n$ and changing $u_n$ by a subsequence if necessary we may assume $u_n\to u$. The geodesic $\alpha(t)=exp_qtu$, $t\in\R$, is a ray. Then $dK_q(u)=a_\alpha\leq 0$ by Corollary 2.2. For each $n$, $K_{\gamma_n}$ is increasing since is an affine function and 
$$
K_{\gamma_n}(0)=K(q)=K(\gamma(t_o))<K(\gamma(t_o+t_n))=K_{\gamma_n}(d(q,\gamma(t_o+t_n)).
$$
Recall that $K_\gamma$ is increasing. Thus $K^\prime_{\gamma_n}(t)=K^\prime_{\gamma_n}(0)=dK_{q}(u_n)> 0$. Taking the limit we deduce that $dK_q(u)=\underset {n\to \infty} \to\lim dK_q(u_n)\geq 0$. Then, $dK_q(u)=0$ 
and $K_\alpha$ is constant non-positive along $\alpha$ by Corollary 2.2. But $K_\alpha(0)=K(q)>0$ give us a contradiction. Hence $K_\gamma$ is constant for any normalized geodesic $\gamma$. Now follows that $M^2$ has constant Gaussian curvature since it is complete.\qed  
\enddemo

\enddocument